\documentclass[a4paper,reqno]{amsart}

\usepackage{amssymb}
\usepackage{amstext}
\usepackage{amsmath}
\usepackage{amscd}
\usepackage{amsthm}
\usepackage{amsfonts}
\usepackage{enumerate}
\usepackage{graphicx}
\usepackage{latexsym}
\usepackage[all]{xy}
\input xy
\xyoption{all}
\usepackage[usenames]{color}

\newtheorem{theorem}{Theorem}[section]

\newtheorem{corollary}[theorem]{Corollary}
\newtheorem{lemma}[theorem]{Lemma}
\newtheorem{proposition}[theorem]{Proposition}

\theoremstyle{definition}
\newtheorem{definition}[theorem]{Definition}
\newtheorem{remark}[theorem]{Remark}
\newtheorem{example}[theorem]{Example}

\newcommand{\DD}{{\mathcal D}}
\newcommand{\UU}{{\mathcal U}}
\newcommand{\SSS}{\mathop{\nu}\nolimits}
\newcommand{\Z}{\mathbf{Z}}
\newcommand{\add}{\operatorname{add}\nolimits}
\newcommand{\proj}{\operatorname{proj}\nolimits}
\newcommand{\inj}{\operatorname{inj}\nolimits}
\newcommand{\gl}{\operatorname{gl.dim}\nolimits}
\newcommand{\dom}{\operatorname{dom.dim}\nolimits}
\newcommand{\id}{\operatorname{id}\nolimits}
\newcommand{\Ker}{\operatorname{Ker}\nolimits}

\newcommand{\Hom}{\operatorname{Hom}\nolimits}
\newcommand{\Endm}{\operatorname{End}\nolimits}
\newcommand{\Ext}{\operatorname{Ext}\nolimits}
\newcommand{\Tor}{\operatorname{Tor}\nolimits}
\newcommand{\RHom}{\mathbf{R}\strut\kern-.2em\operatorname{Hom}\nolimits}
\newcommand{\Lotimes}{\mathop{\stackrel{\mathbf{L}}{\otimes}}\nolimits}
\renewcommand{\mod}{\operatorname{mod}\nolimits}

\begin{document}
\title{$n$-representation-finite algebras and twisted fractionally Calabi-Yau algebras}
\author{Martin Herschend and Osamu Iyama}
\address{O. Iyama: Graduate School of Mathematics, Nagoya University, Chikusa-ku, Nagoya,
464-8602 Japan}
\email{iyama@math.nagoya-u.ac.jp}
\address{M. Herschend: Graduate School of Mathematics, Nagoya University, Chikusa-ku, Nagoya,
464-8602 Japan}
\email{martin.herschend@gmail.com}

\begin{abstract}
In this paper, we study $n$-representation-finite algebras from the viewpoint of the fractionally Calabi-Yau
property. We shall show that all $n$-representation-finite algebras
are twisted fractionally Calabi-Yau. 
We also show that for any $\ell>0$, twisted $\frac{n(\ell-1)}{\ell}$-Calabi-Yau
algebras of global dimension at most $n$ are $n$-representation-finite. 
As an application, we give a construction of $n$-representation-finite
algebras using the tensor product. 
\end{abstract}
\maketitle

The Calabi-Yau (CY) property of triangulated categories was introduced by Kontsevich \cite{Ko} (see also \cite{Ke1}).
It has played important roles in the representation theory of algebras, especially in the categorification program of Fomin-Zelevinsky cluster algebras by cluster tilting theory
(e.g. \cite{Am,BIRS,BMRRT,GLS1,GLS2,IR,Ke2,Ke3,KR}).
The derived categories of finite dimensional non-semisimple algebras are never CY, but they are often \emph{fractionally} CY.
Some recent meetings \cite{Ba,Bi,P,T} as well as recent results including \cite{D,DL,IIKNS,IKM,KST,Ke2,Ladk,Le,Mi,MY,Y} suggest
that the fractionally CY property becomes more and more important in representation theory, singularity theory, commutative and non-commutative algebraic geometry.
The aim of this short paper is to apply the fractionally CY property in the study of $n$-representation-finite algebras defined below.

\medskip
There is a natural generalization of the classical notion of representation-finiteness from the viewpoint of a higher analogue of Auslander-Reiten theory \cite{Au,I1,I3}
studied by several authors \cite{EH,GLS1,IO1,IO2,Lada,HZ1,HZ2,HZ3}.

\begin{definition}
We say that a finite dimensional algebra $\Lambda$ over a field $K$ is \emph{$n$-representation-finite} (for a positive integer $n$)
if $\gl\Lambda\le n$ and there exists an \emph{$n$-cluster tilting $\Lambda$-module $M$}, i.e.
\begin{eqnarray*}
\add M&=&\{X\in\mod\Lambda\ |\ \Ext^i_\Lambda(M,X)=0\ \mbox{ for any }\ 0<i<n\},\\
&=&\{X\in\mod\Lambda\ |\ \Ext^i_\Lambda(X,M)=0\ \mbox{ for any }\ 0<i<n\}.
\end{eqnarray*}
\end{definition}

Clearly $1$-cluster tilting $\Lambda$-modules are additive generators of $\mod\Lambda$, and so $\Lambda$ is $1$-representation-finite if and only if it is representation-finite and hereditary.
The following result \cite[1.3(b)]{I4} plays an important role in this paper.

\begin{proposition}\label{tau n}
Let $\Lambda$ be an $n$-representation-finite algebra. 
Let 
\begin{eqnarray*}
\tau_n:=\Tor^\Lambda_n(D\Lambda,-)\simeq D\Ext^n_\Lambda(-,\Lambda)&:&\mod\Lambda\to\mod\Lambda,\\
\tau_n^-:=D\Tor^\Lambda_n(D-,D\Lambda)\simeq \Ext^n_\Lambda(D\Lambda,-)&:&mod\Lambda\to\mod\Lambda.
\end{eqnarray*}
Let $P_1,\ldots,P_a$ be the isoclasses of indecomposable projective $\Lambda$-modules, and let $I_i$ be the indecomposable injective $\Lambda$-module corresponding to $P_i$.
\begin{itemize}
\item[(a)] There exists a permutation $\sigma\in\mathfrak{S}_a$ and positive integers $\ell_1,\ldots,\ell_a$ such that $\tau_n^{\ell_i-1}I_i\simeq P_{\sigma(i)}$ for any $i$.
\item[(b)] There exists a unique basic $n$-cluster tilting $\Lambda$-module $M$, which is given as the direct sum of the following mutually non-isomorphic indecomposable $\Lambda$-modules.
\[\begin{array}{cccccc}
I_1,&\tau_nI_1,&\tau_n^2I_1,&\cdots&\tau_n^{\ell_1-2}I_1,&\tau_n^{\ell_1-1}I_1\simeq P_{\sigma(1)}\\
I_2,&\tau_nI_2,&\tau_n^2I_2,&\cdots&\tau_n^{\ell_2-2}I_2,&\tau_n^{\ell_2-1}I_2\simeq P_{\sigma(2)}\\
\cdots&\cdots&\cdots&\cdots&\cdots&\cdots\\
I_a,&\tau_nI_a,&\tau_n^2I_a,&\cdots&\tau_n^{\ell_a-2}I_a,&\tau_n^{\ell_a-1}I_a\simeq P_{\sigma(a)}
\end{array}\]
\item[(c)] We have mutually quasi-inverse equivalences $\tau_n:\add(M/\Lambda)\simeq\add(M/D\Lambda)$
and $\tau_n^-:\add(M/D\Lambda)\simeq\add(M/\Lambda)$.
\end{itemize}
\end{proposition}


\medskip
Throughout this paper, let $K$ be a field and $\Lambda$ a finite dimensional $K$-algebra of finite global dimension.
All modules are finitely generated left modules unless stated otherwise.
The composition $fg$ of morphisms (or arrows) means
that $f$ is first and $g$ is next.
We denote by $\DD=\DD^{\rm b}(\mod\Lambda)$ the bounded derived category of the category $\mod\Lambda$
of all $\Lambda$-modules. Then the Nakayama functor
\[\SSS=\SSS_\Lambda:=(D\Lambda)\Lotimes_\Lambda-\simeq D\circ\RHom_\Lambda(-,\Lambda):\DD\to\DD\]
gives a Serre functor, i.e. there exists a functorial isomorphism
\[\Hom_{\DD}(X,Y)\simeq D\Hom_{\DD}(Y,\SSS X)\]
for any $X,Y\in\DD$ \cite{H,BK}.
For each $K$-algebra endomorphism $\phi$ of $\Lambda$, we define an endofunctor of $\DD$ by
\[\phi^*:=\Lambda_\phi\Lotimes_\Lambda-:\DD\to\DD,\]
where we denote by $\Lambda_\phi$ the $\Lambda\otimes_K\Lambda^{\rm op}$-module such that the right action is given by $a\cdot b:=a\phi(b)$.
For any $x\in\Ker\phi$, the morphism $\Lambda\to\Lambda$, $\lambda\mapsto \lambda x$ is mapped to $0$ by $\phi^*$.
Thus $\phi^*$ is an autofunctor if and only if $\phi$ is an automorphism.

\begin{definition}
We say that $\Lambda$ is \emph{twisted fractionally CY}
(or \emph{twisted $\frac{m}{\ell}$-CY})
if there exists an isomorphism
\begin{equation}\label{twisted CY}
\SSS^\ell\simeq[m]\circ\phi^*
\end{equation}
of functors for some integers $\ell\neq0$ and $m$ and
a $K$-algebra endomorphism $\phi$ of $\Lambda$.
Then $\phi$ must be an automorphism since $\SSS$ is an autofunctor.
When $\phi=\id$, we say that $\Lambda$ is \emph{fractionally CY}
(or \emph{$\frac{m}{\ell}$-CY}).
\end{definition}

Notice that (twisted) $\frac{m}{\ell}$-CY algebras are also (twisted) $\frac{mk}{\ell k}$-CY for any $k\neq0$.
The converse does not hold in general, but the Calabi-Yau dimension is unique as a rational number.

\medskip
\noindent{\bf Acknowledgement }
The authors would like to thank an anonymous referee for careful reading the first draft and helpful comments.
The second author was supported by JSPS Grant-in-Aid for Scientific Research 21740010 and 60015849.
The first author is grateful to JSPS for funding his stay at Nagoya University, during which this paper was written.

\section{Main results}

Our first result is the following.

\begin{theorem}\label{result1}
Let $\Lambda$ be a ring-indecomposable $n$-representation-finite algebra. Then:
\begin{itemize}
\item[(a)] $\Lambda$ is twisted fractionally CY.
\item[(b)] More precisely, let $a$ be the number of simple $\Lambda$-modules, and let $b$ be 
the number of idecomposable summands of the basic $n$-cluster tilting $\Lambda$-module.
Then $\Lambda$ is twisted $\frac{n(b-a)k}{bk}$-CY for some $k>0$.
\end{itemize}
\end{theorem}

Our second result shows that a certain converse of the statement (a) above holds.

\begin{definition}
We say that an $n$-representation-finite algebra $\Lambda$ is \emph{$\ell$-homogeneous}
(or simply \emph{homogeneous})
if $\ell:=\ell_1=\cdots=\ell_a$ holds in the notation of Proposition \ref{tau n}.
In this case, we have $\frac{b}{a}=\ell$ for $a$ and $b$ in Theorem \ref{result1}(b).
\end{definition}

\begin{theorem}\label{result2}
For a finite dimensional algebra $\Lambda$ and a positive integer $\ell$, the following conditions are equivalent.
\begin{itemize}
\item[(a)] $\Lambda$ is $\ell$-homogeneous $n$-representation-finite.
\item[(b)] $\Lambda$ is twisted $\frac{n(\ell-1)}{\ell}$-CY and $\gl\Lambda\le n$.
\end{itemize}
\end{theorem}

A typical example is given by tubular algebras of type $(2,2,2,2)$ \cite{R} for the case $n=\ell=2$.
They are known to be $\frac{2}{2}$-CY \cite{LM}, and also to be $2$-representation-finite
\cite{IO2}.

\medskip
Let us apply these results to construct $n$-representation-finite algebras using the tensor product.
The next observation shows that fractionally CY algebras behave nicely under the tensor product,
where we denote $\otimes_K$ by $\otimes$.
The requirement on the field $K$ to be perfect is needed to guarantee
that the tensor product of two finite dimensional $K$-algebras of finite global
dimension is still of finite global dimension.

\begin{proposition}\label{tensor CY}
Let $K$ be a perfect field.
If $\Lambda_i$ is $\frac{m_i}{\ell_i}$-CY (respectively, twisted $\frac{m_i}{\ell_i}$-CY),
then $\Lambda_1\otimes\cdots\otimes\Lambda_k$ is $\frac{m}{\ell}$-CY (respectively, twisted $\frac{m}{\ell}$-CY)
for the least common multiple $\ell$ of $\ell_1,\ldots,\ell_k$ and $m:=\ell(\frac{m_1}{\ell_1}+\cdots+\frac{m_k}{\ell_k})$.
\end{proposition}

As an easy consequence, we have the following useful result.

\begin{corollary}\label{result3}
Let $K$ be a perfect field and $\ell$ a positive integer.
If $\Lambda_i$ is $\ell$-homogeneous $n_i$-representation-finite,
then $\Lambda_1\otimes\cdots\otimes\Lambda_k$ is an $\ell$-homogeneous $(n_1+\cdots+n_k)$-representation-finite
algebra with an $(n_1+\cdots+n_k)$-cluster tilting module
$\bigoplus_{i=0}^{\ell-1}(\tau_{n_1}^i\Lambda_1\otimes\cdots\otimes\tau_{n_k}^i\Lambda_k)$.
\end{corollary}

Notice that the statement is not true if we drop the assmption of $\ell$-homogeneity (see Remark \ref{remarks}(a)).

\medskip
The relationship between $n$-representation-finiteness and the factionally CY property is summarized in the following diagram:

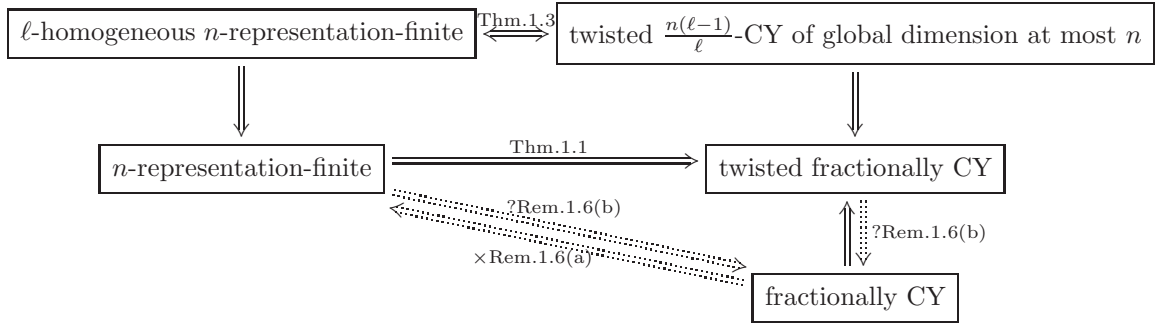
\[\xymatrix{
\fboxsep=1.5\fboxsep\fbox{$\ell$-homogeneous $n$-representation-finite}\ar@{<=>}[r]^(0.447){\rm Thm.1.3}\ar@{=>}[d]&
\fboxsep=1.5\fboxsep\fbox{twisted $\frac{n(\ell-1)}{\ell}$-CY of global dimension at most $n$}\ar@{=>}[d]\\
\fboxsep=1.5\fboxsep\fbox{$n$-representation-finite}\ar@<.7ex>@{=>}[r]^{\rm Thm. 1.1}
\ar@<.7ex>@{:>}[rd]^{\rm ? Rem. 1.6(b)}&
\fboxsep=1.5\fboxsep\fbox{twisted fractionally CY}\ar@<.7ex>@{:>}[d]^{\rm ? Rem. 1.6(b)}\\
&\fboxsep=1.5\fboxsep\fbox{fractionally CY}
\ar@<.7ex>@{=>}[u]\ar@<.7ex>@{:>}[lu]^{\rm \times Rem. 1.6(a)}
}
\]

\begin{remark}\label{remarks}
\begin{itemize}
\item[(a)] Let $\Lambda:=K[\bullet\to\bullet]\otimes K[\bullet\to\bullet]$.
Then $\Lambda$ is fractionally CY, but not $n$-representation-finite for any $n$.
Moreover $\Lambda$ is derived equivalent to the path algebras of type $D_4$, which are $1$-representation-finite.
In particular we have the following conclusions.
\begin{itemize}
\item[$\bullet$] There exist fractionally CY algebras which are not $n$-representation-finite for any $n$.
\item[$\bullet$] $n$-representation-finiteness is not preserved under derived equivalence.
\item[$\bullet$] $n$-representation-finiteness is not preserved under the tensor product.
\end{itemize}
\item[(b)] Although we already know that every $n$-representation-finite algebra is twisted fractionally CY,
we do not know the answer to the following question.
\begin{itemize}
\item[$\bullet$] Is every $n$-representation-finite algebra fractionally CY?
\end{itemize}
Also we do not know the answer to the following related question.
\begin{itemize}
\item[$\bullet$] Is every twisted fractionally CY algebra fractionally CY?
\end{itemize}
This is equivalent to that $\phi$ in \eqref{twisted CY} has a finite order in the outer automorphism group of $\Lambda$.
\item[(c)] Since the Serre functor is unique up to isomorphism, the fractionally CY property is invariant
under derived equivalence. However, we do not know the answer to the following question.
\begin{itemize}
\item[$\bullet$] Are twisted fractionally CY algebras closed under derived equvalence?
\end{itemize}
If we fix the CY dimension, then the answer is negative.
For example $K[\bullet\to\bullet\leftarrow\bullet]$ is twisted $\frac{1}{2}$-CY, but
$K[\bullet\to\bullet\to\bullet]$ is not (see section 3.1).
\end{itemize}
\end{remark}

We end this section by an observation on endomorphism algebras $\Gamma:=\Endm_\Lambda(M)$ of $n$-cluster tilting $\Lambda$-modules $M$.
They are called \emph{$n$-Auslander algebras}, and satisfies $\gl\Gamma\le n+1\le\dom\Gamma$ \cite{I2}.
We have the following descrption of $\Gamma$ when $\Lambda$ is $\ell$-homogeneous $n$-representation-finite.

\begin{proposition}\label{description of Auslander}
If $\Lambda$ in an $\ell$-homogeneous $n$-representation-finite algebra, then the $n$-Auslander algebra of $\Lambda$ is described as
\[\def\arraystretch{.5}\left(\begin{array}{ccccc}
\Lambda&T&\cdots&T^{\otimes_{\Lambda}(\ell-2)}&T^{\otimes_{\Lambda}(\ell-1)}\\
0&\Lambda&\cdots&T^{\otimes_{\Lambda}(\ell-3)}&T^{\otimes_{\Lambda}(\ell-2)}\\
\vdots&\vdots&\ddots&\vdots&\vdots\\
0&0&\cdots&\Lambda&T\\
0&0&\cdots&0&\Lambda
\end{array}\right)\]
for the $\Lambda\otimes\Lambda^{\rm op}$-module $T:=\Ext^{n}_{\Lambda}(D\Lambda,\Lambda)$.
\end{proposition}

By Corollary \ref{result3}, the $(n_1+\cdots+n_k)$-Auslander algebra of $\Lambda_1\otimes\cdots\otimes\Lambda_k$ is given by the `component-wise tensor product'
\[\def\arraystretch{.5}\left(\begin{array}{ccccc}
\Lambda_1\otimes\cdots\otimes\Lambda_k&T_1\otimes\cdots\otimes T_k&\cdots&T_1^{\otimes_{\Lambda_1}(\ell-2)}\otimes\cdots\otimes T_k^{\otimes_{\Lambda_k}(\ell-2)}&T_1^{\otimes_{\Lambda_1}(\ell-1)}\otimes\cdots\otimes T_k^{\otimes_{\Lambda_k}(\ell-1)}\\
0&\Lambda_1\otimes\cdots\otimes\Lambda_k&\cdots&T_1^{\otimes_{\Lambda_1}(\ell-3)}\otimes\cdots\otimes T_k^{\otimes_{\Lambda_k}(\ell-3)}&T_1^{\otimes_{\Lambda_1}(\ell-2)}\otimes\cdots\otimes T_k^{\otimes_{\Lambda_k}(\ell-2)}\\
\vdots&\vdots&\ddots&\vdots&\vdots\\
0&0&\cdots&\Lambda_1\otimes\cdots\otimes\Lambda_k&T_1\otimes\cdots\otimes T_k\\
0&0&\cdots&0&\Lambda_1\otimes\cdots\otimes\Lambda_k
\end{array}\right)\]
of the $n_i$-Auslander algebras.
This is interesting since in general the tensor product of higher Auslander algebras is not a higher Auslander algebra again.

\section{Homogeneity and $(n+1)$-preprojective algebras}\label{section2}

In this section, we further study the homogeneity of $n$-representation-finite algebras.

The following result gives a useful criterion for $n$-representation-finite
algebras to be homogeneous.

\begin{proposition}\label{homogeneous}
Let $\Lambda$ be a ring-indecomposable $n$-representation-finite algebra.
Then $\Lambda$ is homogeneous if and only if $\ell_i=\ell_{\sigma(i)}$ holds
for any $i$.
\end{proposition}


We put
\[\SSS_n:=\SSS\circ[-n]:\DD\to\DD.\]
For an arbitrary $n$-representation-finite algebra $\Lambda$, we denote by
\[\Pi:=\bigoplus_{i\in\Z}\Hom_{\DD}(\Lambda,\SSS_n^{-i}\Lambda)\]
the corresponding \emph{$(n+1)$-preprojective algebra}.
This is a positively graded $K$-algebra, and the multiplication of $f\in\Hom_{\DD}(\Lambda,\SSS_n^{-i}\Lambda)$ and
$g\in\Hom_{\DD}(\Lambda,\SSS_n^{-j}\Lambda)$ is given by
\[f\cdot g:=f\SSS_n^{-i}g.\]
Then $\Pi$ can be viewed as the tensor algebra over $\Lambda$ of the bimodule
$\Ext^n_\Lambda(D\Lambda,\Lambda)$ \cite[2.12]{IO1}.
Notice that $\Pi$ is the cohomology in degree zero of the DG algebra defined in \cite{Ke3}.
It is shown in \cite{IO2} that $\Pi$ is a finite dimensional selfinjective algebra.

A \emph{Nakayama automorphism} $\psi$ of $\Pi$ is an automorphism which
gives an isomorphism $D\Pi\simeq\Pi_\psi$ of $\Pi\otimes\Pi^{\rm op}$-modules.
This is uniquely determined as an element of the outer automorphism group of $\Pi$.
Moreover
\[\psi^*=\Pi_\psi\otimes_\Pi-\simeq(D\Pi)\otimes_\Pi-\simeq D\Hom_\Pi(-,\Pi):\mod\Pi\to\mod\Pi\]
gives the Nakayama functor.

The next result gives the relationship between $\psi$ and $\sigma$.
Notice that the isoclasses of indecomposable projective $\Pi$-modules are give by $\Pi\otimes_\Lambda P_i$ for $1\le i\le a$.

\begin{proposition}\label{psi and sigma}
Let $\psi$ be a Nakayama automorphism of $\Pi$, and let $\sigma$ be as in Proposition \ref{tau n}.
Then we have $\psi^*(\Pi\otimes_\Lambda P_i)\simeq\Pi\otimes_\Lambda P_{\sigma(i)}$ for any $i$.
In particular $\sigma$ gives the Nakayama permutation of $\Pi$.
\end{proposition}

The following main result in this section characterizes homogeneity of $n$-representation-finite
algebras in terms of the corresponding $(n+1)$-preprojective algebras.

\begin{theorem}\label{characterization}
Let $\psi$ be a Nakayama automorphism of $\Pi$, and let $I$ be
the ideal of $\Pi$ such that $\Lambda=\Pi/I$.
Then $\Lambda$ is homogeneous if and only if $\psi(I)=I$.
\end{theorem}


\medskip
When $\Lambda$ is homogeneous, we can interpret a Nakayama automorphism of $\Pi$ as follows:
By Theorem \ref{result2}, there exists a $K$-algebra automorphism $\phi$ of $\Lambda$ such
that $\SSS^\ell\simeq[n(\ell-1)]\circ\phi^*$ for some $\ell\ge1$.

\begin{proposition}\label{phi is Nakayama}
Assume that $\Lambda$ is homogeneous. Then the automorphism $\phi$ of $\Lambda$ extends to a graded $K$-algebra automorphism $\widetilde{\phi}$ of $\Pi$
which is a Nakayama automorphism.
\end{proposition}

\section{Examples}

\subsection{Path algebras}
Let us consider Dynkin diagrams:
\[\xymatrix@C0.4cm@R0.2cm{
A_n&1\ar@{-}[r]&2\ar@{-}[r]&3\ar@{-}[r]&\ar@{..}[rrr]&&&\ar@{-}[r]&n-2\ar@{-}[r]&n-1\ar@{-}[r]&n\\
&&&&&&&&n\\
D_n&1\ar@{-}[r]&2\ar@{-}[r]&3\ar@{-}[r]&\ar@{..}[rrr]&&&\ar@{-}[r]&n-2\ar@{-}[r]\ar@{-}[u]&n-1\\
&&&4\\
E_6&1\ar@{-}[r]&2\ar@{-}[r]&3\ar@{-}[r]\ar@{-}[u]&5\ar@{-}[r]&6\\
&&&&7\\
E_7&1\ar@{-}[r]&2\ar@{-}[r]&3\ar@{-}[r]&4\ar@{-}[r]\ar@{-}[u]&5\ar@{-}[r]&6\\
&&&&&8\\
E_8&1\ar@{-}[r]&2\ar@{-}[r]&3\ar@{-}[r]&4\ar@{-}[r]&5\ar@{-}[r]\ar@{-}[u]&6\ar@{-}[r]&7}\]

The Coxeter number $h$ of each Dynkin diagram is given as follows:
\[\begin{array}{|c|c|c|c|c|}
\hline
A_n&D_n&E_6&E_7&E_8\\ \hline
n+1&2(n-1)&12&18&30\\ \hline
\end{array}\]

The following proposition is well-known \cite[4.1]{MY}.

\begin{proposition}
Let $Q$ be an acyclic quiver. Then $KQ$ is (twisted) fractionally CY
if and only if $Q$ is a Dynkin quiver. In this case, $KQ$ is $\frac{h-2}{h}$-CY,
and if $Q$ is $D_n$ with even $n$, $E_7$ or $E_8$,
then $KQ$ is $\frac{\frac{h}{2}-1}{\frac{h}{2}}$-CY.
\end{proposition}

We define an involution $\omega$ of each Dynkin diagram as follows:
\begin{itemize}
\item For $A_n$, we put $\omega(i)=n+1-i$.
\item For $D_n$ with odd $n$, we put $\omega(n-1)=n$, $\omega(n)=n-1$ and $\omega(i)=i$ for other $i$.
\item For $E_6$, we put $\omega(1)=6$, $\omega(2)=5$, $\omega(5)=2$, $\omega(6)=1$ and $\omega(i)=i$ for other $i$.
\item For other types, we put $\omega=\id$.
\end{itemize}

\begin{proposition}\label{1-rep-fin of type ell}
Let $Q$ be an acyclic quiver.
\begin{itemize}
\item[(a)] The following conditions are equivalent.
\begin{itemize}
\item[(i)] $KQ$ is $\ell$-homogeneous $1$-representation-finite for some $\ell$.
\item[(ii)] $Q$ is a Dynkin quiver and the orientation is stable under $\omega$.
\end{itemize}
\item[(b)] If the conditions in (a) are satisfied, then $\ell=\frac{h}{2}$ and $KQ$ is twisted $\frac{\frac{h}{2}-1}{\frac{h}{2}}$-CY
and the twist is induced by $\omega$.
In particular, if $\omega=\id$, then $KQ$ is $\frac{\frac{h}{2}-1}{\frac{h}{2}}$-CY.
\item[(c)] The underlying diagrams of $\omega$-stable Dynkin quivers are classified by $\ell$ as follows.
\[\begin{array}{|c|c|c|c|c|}
\hline                                                                         
\ell& \ell\ (\neq6,9,15) &6 & 9 & 15 \\ \hline
Q &A_{2\ell-1}, D_{\ell+1}& A_{11}, D_7, E_6& A_{17}, D_{10}, E_7&A_{29}, D_{16},E_8\\ \hline
\end{array}\]        
\end{itemize}
\end{proposition}

\begin{proof}
(a) This follows from Gabriel's description of Auslander-Reiten quivers of $KQ$ \cite{G}.

(b) The equality $\ell=\frac{h}{2}$ is well-known \cite{G}, and the latter assertion follows from Theorem \ref{result2}.
\end{proof}

By Corollary \ref{result3}, for any choice of homogeneous algebras $\Lambda_1,\ldots,\Lambda_k$ which
belong to the same column in the table in (c) above, we have an $\ell$-homogeneous
$k$-representation-finite algebra $\Lambda_1\otimes\cdots\otimes\Lambda_k$.

As one of the simplest examples, we have the following $2$-representation-finite algebras.
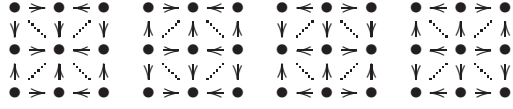
\[\xymatrix@C0.2cm@R0.2cm{
\bullet\ar[r]\ar[d]&\bullet\ar[d]&\bullet\ar[l]\ar[d]&
  \bullet\ar[r]&\bullet\ar@{..}[ld]\ar@{..}[rd]&\bullet\ar[l]&
    \bullet\ar[d]&\bullet\ar[r]\ar[l]\ar[d]&\bullet\ar[d]&
      \bullet\ar@{..}[rd]&\bullet\ar[r]\ar[l]&\bullet\ar@{..}[ld]\\
\bullet\ar[r]&\bullet\ar@{..}[lu]\ar@{..}[ld]\ar@{..}[ru]\ar@{..}[rd]&\bullet\ar[l]&
  \bullet\ar[r]\ar[d]\ar[u]&\bullet\ar[d]\ar[u]&\bullet\ar[l]\ar[d]\ar[u]&
    \bullet\ar@{..}[ru]\ar@{..}[rd]&\bullet\ar[r]\ar[l]&\bullet\ar@{..}[lu]\ar@{..}[ld]&
      \bullet\ar[d]\ar[u]&\bullet\ar[r]\ar[l]\ar[d]\ar[u]&\bullet\ar[d]\ar[u]\\
\bullet\ar[r]\ar[u]&\bullet\ar[u]&\bullet\ar[l]\ar[u]&
  \bullet\ar[r]&\bullet\ar@{..}[lu]\ar@{..}[ru]&\bullet\ar[l]&
    \bullet\ar[u]&\bullet\ar[r]\ar[l]\ar[u]&\bullet\ar[u]&
      \bullet\ar@{..}[ru]&\bullet\ar[r]\ar[l]&\bullet\ar@{..}[lu]}\]

\subsection{$n$-representation-finite algebras of type $A$}

We give a different class of examples of $\ell$-homogeneous $n$-representation-finite algebras
based on results in \cite{IO1}.
Fix positive integers $n$ and $s$. We define the quiver $Q:=Q^{(n,s)}$ 
with the set $Q_0$ of vertices and the set $Q_1$ of arrows by
\begin{eqnarray*}
Q_0&:=&\{x=(x_1,x_2,\ldots,x_{n+1})\in\Z_{\ge0}^{n+1}\ |\ \sum_{i=1}^{n+1}x_i=s-1\},\\
Q_1&:=&\{x\stackrel{i}{\to}x+f_i\ |\ 1\le i\le n+1,\ x,x+f_i\in Q_0\}
\end{eqnarray*}
where $f_i$ denotes the vector $f_i:=(0,\ldots,0,\stackrel{i}{-1},\stackrel{i+1}{1},0,\ldots,0)$ for $1\le i\le n$ and $f_{n+1}:=(\stackrel{1}{1},0,\ldots,0,\stackrel{n+1}{-1})$.

\begin{example}
$Q^{(1,5)}$, $Q^{(2,4)}$ and $Q^{(3,3)}$ are the following quivers.
\[\xymatrix@C0cm@R0.2cm{
{\scriptstyle 40}\ar@<.5ex>[rr]&&{\scriptstyle 31}\ar@<.5ex>[rr]\ar@<.5ex>[ll]&&{\scriptstyle 22}\ar@<.5ex>[rr]\ar@<.5ex>[ll]&&{\scriptstyle 13}\ar@<.5ex>[rr]\ar@<.5ex>[ll]&&{\scriptstyle 04}\ar@<.5ex>[ll]&
  &&&{\scriptstyle 030}\ar[rd]&&&&
    &&{\scriptstyle 0200}\ar[rd]&&\\
&&&&&&&&&
  &&{\scriptstyle 120}\ar[ru]\ar[rd]&&{\scriptstyle 021}\ar[ll]\ar[rd]&&&
    &{\scriptstyle 1100}\ar[ru]\ar[rd]&&{\scriptstyle 0110}\ar[rd]\ar[ldd]&\\
&&&&&&&&&
  &{\scriptstyle 210}\ar[ru]\ar[rd]&&{\scriptstyle 111}\ar[ll]\ar[ru]\ar[rd]&&{\scriptstyle 012}\ar[ll]\ar[rd]&&
    {\scriptstyle 2000}\ar[ru]&&{\scriptstyle 1010}\ar[ru]\ar[ldd]&&{\scriptstyle 0020}\ar[ldd]\\
&&&&&&&&&
  {\scriptstyle 300}\ar[ru]&&{\scriptstyle 201}\ar[ll]\ar[ru]&&{\scriptstyle 102}\ar[ll]\ar[ru]&&{\scriptstyle 003}\ar[ll]&
    &&{\scriptstyle 0101}\ar[rd]\ar[luu]&&\\
&&&&&&&&&
  &&&&&&&
    &{\scriptstyle 1001}\ar[ru]\ar[luu]&&{\scriptstyle 0011}\ar[luu]\ar[ldd]&\\
&&&&&&&&&
  &&&&&&&
    &&&&\\
&&&&&&&&&
  &&&&&&&
    &&{\scriptstyle 0002}\ar[luu]&&
}\]
\end{example}

We define the $K$-algebra $\Gamma=\Gamma^{(n,s)}:=KQ/I$, where
$I$ is the ideal defined by the following relations:

For any $x\in Q_0$ and $1\le i,j\le n+1$ satisfying $x+f_i,x+f_i+f_j\in Q_0$,
\[(x\stackrel{i}{\to}x+f_i\stackrel{j}{\to}x+f_i+f_j)=\left\{
\begin{array}{cl}
(x\stackrel{j}{\to}x+f_j\stackrel{i}{\to}x+f_i+f_j)&\mbox{if}\ x+f_j\in Q_0,\\
0&\mbox{otherwise}.
\end{array}
\right.\]
We call a subset $C$ of $Q_1$ a \emph{cut} if it contains exactly one arrow from each cycle of length $n+1$ in $Q$.
In this case, we define a $K$-algebra by
\[\Lambda_C:=\Gamma/\langle C\rangle,\]
for the ideal $\langle C\rangle$ of $\Gamma$ generated by the arrows in $C$.
A main result in \cite{IO1} is the following.
(Notice that the assumption in \cite{IO1} that $K$ is algebraically closed is not necessary.)

\begin{proposition}
For any cut $C$ of $Q$, the algebra $\Lambda_C$ is an $n$-representation-finite algebra
with the $(n+1)$-preprojective algebra $\Gamma$.
\end{proposition}

It is natural to ask when $\Lambda_C$ is homogeneous.
To answer this question we introduce the automorphism $\omega$ of $Q$ defined by
\[\omega(x_1,x_2,\ldots,x_{n+1}):=(x_{n+1},x_1,\ldots,x_n).\]
We leave the proof of the following statement to the reader.

\begin{theorem}\label{nakayama}
The quiver morphism $\omega$ induces a Nakayama automorphism of $\Gamma$.
\end{theorem}

\begin{proof}
For each $x = (x_1,x_2,\ldots,x_{n+1}) \in \Z^{n+1}$ we adopt the convention $x_{i + n +1} = x_i$. Now let $x \in Q_0$ and $1\leq i \leq n+1$ such that $x +f_{i} \in Q_0$. Then $\omega(x +f_{i}) =  \omega(x) +f_{i+1}\in Q_0$ and
\[
\omega(x \stackrel{i}{\to}x+f_i ) = \omega(x)\stackrel{i+1}{\to} \omega(x)+f_{i+1}.
\]
Now assume $x +f_{i}+f_{j} \in Q_0$. Then 
\[
\omega(x \stackrel{i}{\to}x+f_i \stackrel{j}{\to}x+f_i +f_j) = \omega(x) \stackrel{i+1}{\to}\omega(x)+f_{i+1} \stackrel{j+1}{\to}\omega(x)+f_{i+1} +f_{j+1}.
\]
Since $x+f_j \in Q_0$ if and only if $\omega(x)+f_{j+1} \in Q_0$ it follow that $\omega$ induces an automorphism on $\Gamma$.

The relations defining $I$ are only commutativity and zero relations. Thus the set $B$ of non-zero residue classes $p + I \not = I$ of paths $p$ in $Q$ forms a basis of $\Gamma$. Since $\Gamma$ is finite dimensional, $D\Gamma$ has basis $DB : = \{\bar{b}^* | \bar{b} \in B \}$ dual to $B$. We proceed to show that there is a bijection
\[
\phi : DB \simeq B
\]
that extends to a bimodule isomorphism $D\Gamma \simeq \Gamma_{\omega}$.

Define the quiver $\widehat{Q}$ by
\begin{eqnarray*}
\widehat{Q}_0&:=&\{ x\in\Z^{n+1}\ |\ \sum_{i=1}^{n+1}x_i=s-1\},\\
\widehat{Q}_1&:=&\{x\stackrel{i}{\to}x+f_i\ |\ 1\le i\le n+1\}.
\end{eqnarray*}
Hence $Q$ is the full subquiver of $\widehat{Q}$ with vertecies $\{x \in \widehat{Q}_0 \; | \; x_i \ge 0 \mbox{ for all } 1\leq i \leq n+1\}$. Let $\widehat{I}$ be the ideal in $K\widehat{Q}$ defined by the relations
\[
(x\stackrel{i}{\to}x+f_i\stackrel{j}{\to}x+f_i+f_j)=(x\stackrel{j}{\to}x+f_j\stackrel{i}{\to}x+f_i+f_j)
\]
and set $\widehat{\Gamma} = K\widehat{Q}/\widehat{I}$. We have a canonical surjective ring morphism $\pi : \widehat{\Gamma} \to \Gamma$ with kernel
\[
K = \sum_{z \not \in Q_0} \widehat{\Gamma}e_z\widehat{\Gamma},
\]
where $e_z$ denotes the path of length $0$ at $z$.
Since $\widehat{I}$ is defined by commutativity relations, the residue classes of paths form a basis of $\widehat{\Gamma}$. Moreover, the residue classes of paths that are not in $K$ are mapped bijectively to $B$ by $\pi$.

We define a $\Z^{n+1}$-grading $g$ on $\widehat{\Gamma}$ by $(g(x \stackrel{i}{\to} x +f_i))_j = \delta_{ij}$. It is well-defined since the defining relations of $\widehat{I}$ are homogeneous. For paths $p,p'$ in $\widehat{Q}$ we write $p \equiv p'$ if and only if $p-p' \in \widehat{I}$. Let $p$ be a path from $x$ to $y$ of degree $d$. Then $y - x =\sum_{i}d_if_i$. In fact $ p \equiv p_{x,d,y}$, where
\[
p_{x,d,y} = x \stackrel{1}{\to}\cdots \stackrel{1}{\to} x + d_1f_1\stackrel{2}{\to}\cdots \stackrel{2}{\to} x + d_1f_1 + d_2f_2 
\stackrel{3}{\to} \;\; \cdots \;\; \stackrel{n}{\to}
y - d_{n+1}f_{f+1} \stackrel{n+1}{\to}\cdots \stackrel{n+1}{\to} y.
\]
In particular, $p +\widehat{I}$ is determined by $d$ together with either $x$ or $y$. Moreover, for each path $p'$ from $x$ to $z$, there is a path $q$ from $z$ to $y$ such that $p \equiv p'q$ if and only if $g(p')_i \leq d_i$ for all $1\leq i \leq n+1$ (take $q = p_{z,d-g(p'),y}$).

Observe that $p +\widehat{I} \in K$ if and only if there are paths $p'$ from $x$ to $z$ and $q$ from $z$ to $y$ such that $p \equiv p'q$ and $z \not \in Q_0$. By the above observation this is equivalent to that there is $a \in \Z_{\ge0}^{n+1}$ such that $a_j\leq d_j$ for each $1\leq j \leq n+1$ and $x + \sum_{i}a_if_i \not \in Q_0$. This holds if and only if there is $1\leq j \leq n+1$ such that $x_j < d_j$ (take $a_i = \delta_{ij}d_j$).

We conclude that for every path $p$ in $Q$ from $x$ to $z$ such that $p +I \not = I$ there is a path $q$ in $Q$ from $z$ to some vertex $y$ such that $g(pq) = x$. In particular, $pq +I \not = I$ and
\[
y = x + \sum_{i}x_if_i.
\]
Thus
\[
y_j = x_j - x_j + x_{j-1} = x_{j-1} = \omega(x)_j.
\]
Moreover, $q+I$ is determined by $p$ since $g(q) = x - g(p)$. Similarly, for each path $q$ in $Q$ from $z$ to $\omega(x)$ such that $q \not \in I$ there is a path $p$ in $Q$ from $x$ to $z$, unique up to $\equiv$ such that $g(pq) = x$. We define the bijection $\phi : DB \simeq B$ by $\phi((p+I)^*) = q +I$, where $p$ is a path from $x$ to $z$ and $q$ is a path from $z$ such that $g(pq) = x$. It induces a linear bijection
\[
\phi : D\Gamma \simeq \Gamma_{\omega}.
\]
We proceed to show that $\phi$ is a bimodule isomorphism, thus completing the proof.

First consider the following configuration of paths in $Q$.
\[\xymatrix{
x \ar@{-->}[r]^a &y \ar@{-->}[r]^{p'} &z \ar@{-->}[r]^b &w \ar@{-->}[r]^q &\omega(x) \ar@{-->}[r]^{\omega(a)} &\omega(y)
}\]
Then
\[
g(p'bq\omega(a)) - g(ap'bq) = g(\omega(a)) - g(a) =  \sum_{i} g(a)_if_i = y-x.
\]
In particular, $g(ap'bq) = x$ if and only if $g(p'bq\omega(a)) = y$, or equivalently $\phi((ap'b + I)^*) = q + I$ if and only if $\phi((p' +I)^*) = bq\omega(a) + I$.

Now fix paths $\alpha$, $\beta$, $p$ from $x$ to $w$ and $q$ from $w$ to $\omega(x)$ such that $g(pq) = x$. Set $\bar{b} =(p+I)$ so that $\phi(\bar{b}^*)=q+I$. We claim that there is $p'$ such that $p \equiv \alpha p' \beta$ if and only if $\beta q\omega(\alpha) \not \in I$. If such $p'$ exists then $g(\alpha p' \beta q) = x$ and by our observation $g(p' \beta q \omega(\alpha)) = y$. This implies $\beta q\omega(\alpha) \not \in I$. On the other hand, if $\beta q\omega(\alpha) \not \in I$, then there is $p'$ such that $g(p' \beta q \omega(\alpha)) = y$ and $g(\alpha p' \beta q) = x$ which implies $\alpha p' \beta \equiv p$. 

If there is a $p'$ satisfying $p \equiv \alpha p' \beta$, then by the bimodule structure on $D\Gamma$, we have  $\beta \bar{b}^* \alpha = (p' +I)^*$. Otherwise $\beta \bar{b}^* \alpha = 0$. In either case $\phi(\beta \bar{b}^* \alpha) = \beta q \omega(\alpha)+I = \beta \phi(\bar{b}^*)\alpha$ by the above observations.
\end{proof}

We have the following characterization.

\begin{theorem}\label{type A}
Let $C$ be an admissible set of $Q$.
Then the $n$-representation-finite algebra $\Lambda_C$ is homogeneous if and only if $\omega(C)=C$.
In this case, we have $\ell=\frac{s-1}{n+1}+1$.
\end{theorem}

\begin{proof}
The first assertion follows from Theorems \ref{characterization} and \ref{nakayama}.
To show the second assertion we apply Theorem \ref{result1}(a).
We know that $a$ is equal to the number ${s+n-1\choose n}$ of vertices of the quiver $Q^{(n,s)}$,
and $b$ is equal to the number ${s+n\choose n+1}$ of vertices of the quiver
$Q^{(n+1,s)}$ by \cite[6.12]{I4}. Thus $\ell=\frac{b}{a}=\frac{s+n}{n+1}$.
\end{proof}

For the case $n=1$, this result gives the characterization of type $A_s$ quivers given in Proposition
\ref{1-rep-fin of type ell}.

For the case $n=2$ and $s=4$, there are the following five $2$-homogeneous $2$-representation-finite algebras
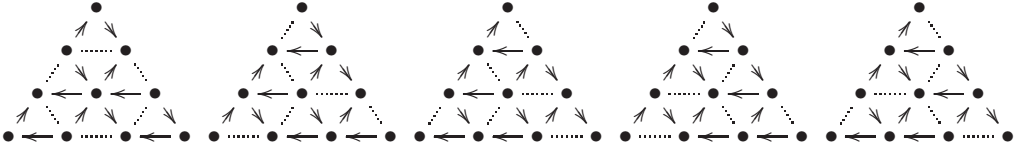
\[\xymatrix@C0cm@R0.2cm{
&&&\bullet\ar[rd]&&&&
  &&&\bullet\ar[rd]&&&&
    &&&\bullet\ar@{..}[rd]&&&&
      &&&\bullet\ar[rd]&&&&
        &&&\bullet\ar@{..}[rd]&&&&\\
&&\bullet\ar[ru]\ar[rd]&&\bullet\ar@{..}[ll]\ar@{..}[rd]&&&
  &&\bullet\ar@{..}[ru]\ar@{..}[rd]&&\bullet\ar[ll]\ar[rd]&&&
    &&\bullet\ar[ru]\ar@{..}[rd]&&\bullet\ar[ll]\ar[rd]&&&
      &&\bullet\ar@{..}[ru]\ar[rd]&&\bullet\ar[ll]\ar[rd]&&&
        &&\bullet\ar[ru]\ar[rd]&&\bullet\ar[ll]\ar[rd]&&&\\
&\bullet\ar@{..}[ru]\ar@{..}[rd]&&\bullet\ar[ll]\ar[ru]\ar[rd]&&\bullet\ar[ll]\ar[rd]&&
  &\bullet\ar[ru]\ar[rd]&&\bullet\ar[ll]\ar[ru]\ar[rd]&&\bullet\ar@{..}[ll]\ar@{..}[rd]&&
    &\bullet\ar[ru]\ar[rd]&&\bullet\ar[ll]\ar[ru]\ar[rd]&&\bullet\ar@{..}[ll]\ar[rd]&&
      &\bullet\ar[ru]\ar[rd]&&\bullet\ar@{..}[ll]\ar@{..}[ru]\ar@{..}[rd]&&\bullet\ar[ll]\ar@{..}[rd]&&
        &\bullet\ar[ru]\ar[rd]&&\bullet\ar@{..}[ll]\ar@{..}[ru]\ar@{..}[rd]&&\bullet\ar[ll]\ar[rd]&&\\
\bullet\ar[ru]&&\bullet\ar[ll]\ar[ru]&&\bullet\ar@{..}[ll]\ar@{..}[ru]&&\bullet\ar[ll]&
  \bullet\ar[ru]&&\bullet\ar@{..}[ll]\ar@{..}[ru]&&\bullet\ar[ll]\ar[ru]&&\bullet\ar[ll]&
    \bullet\ar@{..}[ru]&&\bullet\ar[ll]\ar@{..}[ru]&&\bullet\ar[ll]\ar[ru]&&\bullet\ar@{..}[ll]&
      \bullet\ar[ru]&&\bullet\ar@{..}[ll]\ar[ru]&&\bullet\ar[ll]\ar[ru]&&\bullet\ar[ll]&
        \bullet\ar@{..}[ru]&&\bullet\ar[ll]\ar[ru]&&\bullet\ar[ll]\ar[ru]&&\bullet\ar@{..}[ll]&
}\]

\section{Proof of our results}

Recall that $\Lambda$ is a finite dimensonal $K$-algebra of finite global dimension and
$\DD=\DD^{\rm b}(\mod\Lambda)$ is the bounded derived category of $\mod\Lambda$.

We need the following observation.

\begin{lemma}\label{module version}
For a $\Lambda\otimes\Lambda^{\rm op}$-module $X$, the following conditions are equivalent.
\begin{itemize}
\item[(a)] $X\simeq\Lambda$ as $\Lambda$-modules.
\item[(b)] $X\simeq\Lambda_\phi$ as $\Lambda\otimes\Lambda^{\rm op}$-modules for some $K$-algebra endomorphism $\phi$ of $\Lambda$.
\end{itemize}
\end{lemma}

\begin{proof}
(b)$\Rightarrow$(a) This is clear.

(a)$\Rightarrow$(b) Fix an isomorphism $f\in\Hom_\Lambda(X,\Lambda)$.
For each $\lambda\in\Lambda$, the right multiplication $(\cdot\lambda)\in\Endm_\Lambda(X)$
gives an endomorphism $f^{-1}\cdot(\cdot\lambda)\cdot f\in\Endm_\Lambda(\Lambda)$.
This is given by the right multiplication of an element $\phi(\lambda)\in\Lambda$.
Then $\phi$ gives the desired $K$-algebra endomorphism of $\Lambda$.
\end{proof}

We need the following observation.

\begin{lemma}\label{isomorphisms}
Let $X\in\DD^{\rm b}(\mod\Lambda\otimes\Lambda^{\rm op})$.
\begin{itemize}
\item[(a)] The following conditions are equivalent.
\begin{itemize}
\item[(i)] $X\simeq\Lambda$ in $\DD^{\rm b}(\mod\Lambda)$.
\item[(ii)] $X\simeq\Lambda_\phi$ in $\DD^{\rm b}(\mod\Lambda\otimes\Lambda^{\rm op})$
for some $K$-algebra endomorphism $\phi$ of $\Lambda$.
\end{itemize}
\item[(b)] The following conditions are equivalent.
\begin{itemize}
\item[(i)] There exists an isomorphism $f:\Lambda\to X$ in $\DD^{\rm b}(\mod\Lambda)$
which commutes with the right action of $\Lambda$.
\item[(ii)] $X\simeq\Lambda$ in $\DD^{\rm b}(\mod\Lambda\otimes\Lambda^{\rm op})$.
\end{itemize}
\end{itemize}
\end{lemma}

\begin{proof}
(a) The condition (i) is equivalent to the following one:
\begin{itemize}
\item $H^i(X)=0$ for any $i\neq0$ and $H^0(X)\simeq\Lambda$ as $\Lambda$-modules.
\end{itemize}
In view of Lemma \ref{module version}, this is equivalent to the following condition:
\begin{itemize}
\item $H^i(X)=0$ for any $i\neq0$ and $H^0(X)\simeq\Lambda_\phi$ as $\Lambda\otimes\Lambda^{\rm op}$-modules
for some $K$-algebra endomorphism $\phi$ of $\Lambda$.
\end{itemize}
This is equivalent to the condition (ii).

(b) The condition (i) is equivalent to the following one:
\begin{itemize}
\item $H^i(X)=0$ for any $i\neq0$ and $H^0(X)\simeq\Lambda$ as $\Lambda\otimes\Lambda^{\rm op}$-modules.
\end{itemize}
This is equivalent to the condition (ii).
\end{proof}

Let us show the following characterization of (twisted) fractionally CY property.

\begin{proposition}\label{one side to two side}
Let $\Lambda$ be a finite dimensional $K$-algebra of finite global dimension.
\begin{itemize}
\item[(a)] The following conditions are equivalent.
\begin{itemize}
\item[(i)] $\Lambda$ is twisted $\frac{m}{\ell}$-CY.
\item[(ii)] $\SSS^\ell\Lambda\simeq\Lambda[m]$ in $\DD^{\rm b}(\mod\Lambda)$.
\item[(iii)] There exists an isomorphism $(D\Lambda)^{\Lotimes_\Lambda\ell}\simeq\Lambda[m]$ in $\DD^{\rm b}(\mod\Lambda)$.
\end{itemize}
\item[(b)] The following conditions are equivalent.
\begin{itemize}
\item[(i)] $\Lambda$ is $\frac{m}{\ell}$-CY.
\item[(ii)] $\SSS^\ell\simeq[m]$ as functors on $\DD^{\rm b}(\mod\Lambda)$.
\item[(iii)] There exists an isomorphism $(D\Lambda)^{\Lotimes_\Lambda\ell}\simeq\Lambda[m]$ in $\DD^{\rm b}(\mod\Lambda\otimes\Lambda^{\rm op})$.
\end{itemize}
\end{itemize}
\end{proposition}

\begin{proof}
(a) By definition (ii) is equivalent to (iii). Since $\Lambda_{\phi} \Lotimes_\Lambda \Lambda [m] \simeq \Lambda[m]$, (i) implies (ii).

If (iii) holds, then $(D\Lambda)^{\Lotimes_\Lambda\ell}\simeq\Lambda_\phi[m]$ in $\DD^{\rm b}(\mod\Lambda\otimes\Lambda^{\rm op})$
for some $K$-algebra endomorphism $\phi$ of $\Lambda$ by Lemma \ref{isomorphisms}(a).
Thus we have $\SSS^\ell=[m]\circ\phi^*$, and (i) holds.

(b) By definition (i) is equivalent to (ii). Since $\SSS^\ell=(D\Lambda)^{\Lotimes_\Lambda\ell}\Lotimes_\Lambda-$ and $[m] \simeq \Lambda[m] \Lotimes_\Lambda- $, (iii) implies (ii).
Finally, if (ii) holds, then there exists an isomorphism $f:(D\Lambda)^{\Lotimes_\Lambda\ell}\to\Lambda[m]$ in $\DD^{\rm b}(\mod\Lambda)$
which commutes with the right action of $\Lambda$.
By Lemma \ref{isomorphisms}(b), we have that (iii) holds.
\end{proof}

%

\medskip
We need the following observations.

\begin{lemma}\label{tau and S}
\begin{itemize}
\item[(a)] Let $\Lambda$ be a finite dimensional $K$-algebra with $\gl\Lambda\le n$ and $X\in\mod\Lambda$. Then:
\begin{itemize}
\item[$\bullet$] $\tau_nX\simeq\SSS_nX$ if and only if $\Ext^i_\Lambda(X,\Lambda)=0$ for any $i\neq n$.
\item[$\bullet$] If $X$ is indecomposable non-projective and satisfies $\Ext^i_\Lambda(X,\Lambda)=0$ for any $0<i<n$, then we have $\Hom_\Lambda(X,\Lambda)=0$ and $\tau_nX\simeq\SSS_nX$.
\item[$\bullet$] $\tau_n^-X\simeq\SSS_n^{-1}X$ if and only if $\Ext^i_\Lambda(D\Lambda,X)=0$ for any $i\neq n$.
\item[$\bullet$] If $X$ is indecomposable non-injective and satisfies $\Ext^i_\Lambda(D\Lambda,X)=0$ for any $0<i<n$, then we have $\Hom_\Lambda(D\Lambda,X)=0$ and $\tau_n^-X\simeq\SSS_n^{-1}X$.
\end{itemize}
\item[(b)] Let $\Lambda$ be an $n$-representation-finite algebra and
$M$ an $n$-cluster tilting $\Lambda$-module. Then we have 
\[\tau_nX\simeq\SSS_nX\ \mbox{ and }\ \Hom_\Lambda(X,\Lambda)=0\ \ \ (\mbox{respectively, }\ \tau_n^-X\simeq\SSS_n^{-1}X\ \mbox{ and }\ \Hom_\Lambda(D\Lambda,X)=0)\]
for any $X\in\add M$ without non-zero projective (respectively, injective) direct summands.
\item[(c)] Let $\Lambda$ be a finite dimensional $K$-algebra.
For any $K$-algebra automorphism $\phi$ of $\Lambda$,
we have an isomorphism $\phi^*\SSS\simeq\SSS\phi^*$ of autofunctors of $\DD$.
\end{itemize}
\end{lemma}

\begin{proof}
(a) The first and the third assertions are clear. For the remaining assertions, see \cite[2.3(b)]{I4}.

(b) This is immediate from (a) and the definition of $n$-cluster tilting modules.

(c) There exist isomorphisms
\[\Lambda_\phi\otimes_\Lambda(D\Lambda)\simeq
{}_{\phi^{-1}}(D\Lambda)\simeq D(\Lambda_{\phi^{-1}})\simeq
D({}_\phi\Lambda)\simeq(D\Lambda)_\phi\simeq
(D\Lambda)\otimes_\Lambda\Lambda_\phi\]
of $\Lambda\otimes\Lambda^{\rm op}$-modules.
This gives an isomorphism $\Lambda_\phi\Lotimes_\Lambda(D\Lambda)\simeq(D\Lambda)\Lotimes_\Lambda\Lambda_\phi$
in $\DD^{\rm b}(\mod\Lambda\otimes\Lambda^{\rm op})$.
Thus we have the desired isomorphism.
\end{proof}

\subsection{Proof of Theorem \ref{result1}}

We use the notation in Proposition \ref{tau n}.
By Lemma \ref{tau and S}(b), we have
\[\SSS_n^{\ell_i-1}I_i\simeq P_{\sigma(i)}\]
for any $1\le i\le a$. In particular, we have
\[\SSS^{\ell_i}P_i\simeq\SSS^{\ell_i-1}I_i\simeq P_{\sigma(i)}[n(\ell_i-1)]\]
for any $1\le i\le a$.
We take $c>0$ such that $\sigma^c=\id$, and let
\[m_i:=\ell_i+\ell_{\sigma(i)}+\cdots+\ell_{\sigma^{c-1}(i)}.\]
Then one can easily check that
\[\SSS^{m_i}P_i\simeq P_i[n(m_i-c)]\]
for any $1\le i\le a$.

(a) By Proposition \ref{one side to two side}(a), we only have to check $m_1=\cdots=m_a$.
Since $\Lambda$ is ring-indecomposable, it is enough to check that
$m_i=m_j$ holds if $\Hom_\Lambda(P_i,P_j)\neq0$.
In this case, we have 
\[0\neq\Hom_{\DD}(\SSS^{m_im_j}P_i,\SSS^{m_im_j}P_j)
=\Hom_{\DD}(P_i[n(m_i-c)m_j],P_j[n(m_j-c)m_i])
=\Hom_{\DD}(P_i,P_j[nc(m_j-m_i)]).\]
Thus we have $nc(m_j-m_i)=0$ and $m_i=m_j$.

(b) We have shown that $\Lambda$ is twisted $\frac{n(m-c)}{m}$-CY
for $m:=m_1=\cdots=m_a$. Since
\[\frac{c}{m}=\frac{ac}{m_1+\cdots+m_a}=\frac{a}{\ell_1+\cdots+\ell_a}=\frac{a}{b},\]
we have the assertion.
\qed

\subsection{Proof of Theorem \ref{result2}}

(a)$\Rightarrow$(b) By Lemma \ref{tau and S}(b), we have
\[\SSS_n^{\ell-1}(D\Lambda)\simeq\tau_n^{\ell-1}(D\Lambda)\simeq\Lambda.\]
Thus we have
\[\SSS^\ell\Lambda\simeq\SSS^{\ell-1}(D\Lambda)\simeq\Lambda[n(\ell-1)].\]
The assertion follows from Proposition \ref{one side to two side}(a).

(b)$\Rightarrow$(a) 
We need the following results \cite[3.1]{IO2} and \cite[5.4(a)]{I4}.

\begin{proposition}\label{selfinjective}
Let $\Lambda$ be a finite dimensional $K$-algebra with $\gl\Lambda\le n$.
\begin{itemize}
\item[(a)] Let $\UU:=\add\{\SSS_n^j\Lambda\ |\ j\in\Z\}\subset\DD$.
Then $\Lambda$ is $n$-representation-finite if and only if $D\Lambda\in\UU$.
\item[(b)] Let $(\DD^{\le0},\DD^{\ge0})$ be the standard t-structure of $\DD$.
Then $\SSS_n\DD^{\ge0}\subset\DD^{\ge0}$ and $\SSS_n^{-1}\DD^{\le0}\subset\DD^{\le0}$.
\end{itemize}
\end{proposition}

%

Since $\SSS^\ell\Lambda\simeq \Lambda[n(\ell-1)]$, we have
$\SSS_n^{\ell-1}(D\Lambda)\simeq\Lambda$, which implies $D\Lambda\in\UU$.
By Proposition \ref{selfinjective}(a), we have that $\Lambda$ is $n$-representation-finite.
On the other hand, we have
\[\SSS_n^{-i}\Lambda\simeq\SSS_n^{\ell-1-i}(D\Lambda)\in\DD^{\le0}\cap\DD^{\ge0}=\mod\Lambda\]
for any $0\le i<\ell$ by Proposition \ref{selfinjective}(b).
By Lemma \ref{tau and S}(a), we have $\tau_n^{\ell-1}(D\Lambda)\simeq
\SSS_n^{\ell-1}(D\Lambda)\simeq\Lambda$. Thus $\Lambda$ is $\ell$-homogeneous.
\qed

\subsection{Proof of Proposition \ref{tensor CY} and Corollary \ref{result3}}

We shall prove Proposition \ref{tensor CY}.

Since $K$ is a perfect field, $\Lambda_1\otimes\cdots\otimes\Lambda_k$ has finite global dimension.
We prove the assertion for twisted fractionally CY property.
One can show the assertion for fractionally CY property similarly by
using Proposition \ref{one side to two side}(b) instead of (a).

By Proposition \ref{one side to two side}(a), we have an isomorphism
\[(D\Lambda_i)^{\Lotimes_{\Lambda_i}\ell}\simeq\Lambda_i[n_i(\ell-1)]\]
in $\DD^{\rm b}(\mod\Lambda_i)$.
Thus we have
\begin{eqnarray*}
&&(D(\Lambda_1\otimes\cdots\otimes\Lambda_k))^{\Lotimes_{\Lambda_1\otimes\cdots\otimes\Lambda_k}\ell}
\simeq(D\Lambda_1)^{\Lotimes_{\Lambda_1}\ell}\otimes\cdots\otimes(D\Lambda_k)^{\Lotimes_{\Lambda_k}\ell}\\
&\simeq&\Lambda_1\biggl[\frac{\ell m_1}{\ell_1}\biggr]\otimes\cdots\otimes\Lambda_k\biggl[\frac{\ell m_k}{\ell_k}\biggr]
\simeq(\Lambda_1\otimes\cdots\otimes\Lambda_k)\biggl[\ell\biggl(\frac{m_1}{\ell_1}+\cdots\frac{m_k}{\ell_k}\biggr)\biggr]
\end{eqnarray*}
in $\DD^{\rm b}(\mod\Lambda_1\otimes\cdots\otimes\Lambda_k)$.
By Proposition \ref{one side to two side}(a), we have the assertion.
\qed

\medskip
We shall prove Corollary \ref{result3}.

First of all, note that $\Lambda_1\otimes\cdots\otimes\Lambda_k$ has global dimension at most $n_1+\cdots+n_k$ since $K$ is a perfect field.
By Theorem \ref{result2}(a)$\Rightarrow$(b), we have that
$\Lambda_i$ is twisted $\frac{n_i(\ell-1)}{\ell}$-CY.
By Proposition \ref{tensor CY}, we have that $\Lambda_1\otimes\cdots\otimes\Lambda_k$
is twisted $\frac{(n_1+\cdots+n_k)(\ell-1)}{\ell}$-CY.
By Theorem \ref{result2}(b)$\Rightarrow$(a), we have that
$\Lambda_1\otimes\cdots\otimes\Lambda_k$ is $\ell$-homogeneous $(n_1+\cdots+n_k)$-representation-finite. 

The second statement is clear from Proposition \ref{tau n} since we have
\[\tau_{n_1+\cdots+n_k}(X_1\otimes\cdots\otimes X_k)=(\tau_{n_1}X_1)\otimes\cdots\otimes(\tau_{n_k}X_k)\]
for any $X_i\in\mod\Lambda_i$.
\qed

\medskip
We end this section by giving another proof of Corollary \ref{result3}.

Since $\Lambda_i$ is $\ell$-homogeneous $n_i$-representation-finite, we have
$D\Lambda_i\simeq\nu_{n_i}^{-\ell}\Lambda_i$ for any $i$. Thus we have 
\[D(\Lambda_1\otimes\cdots\otimes\Lambda_k)
\simeq D\Lambda_1\otimes\cdots\otimes D\Lambda_k
\simeq \nu_{n_1}^{-\ell}\Lambda_1\otimes\cdots\otimes\nu_{n_k}^{-\ell}\Lambda_k
\simeq \nu_{n_1+\cdots+n_k}^{-\ell}(\Lambda_1\otimes\cdots\otimes\Lambda_k).\]
Thus $\Lambda_1\otimes\cdots\otimes\Lambda_k$ is $(n_1+\cdots+n_k)$-representation-finite
by Proposition \ref{selfinjective}(a).
\qed

\subsection{Proof of Proposition \ref{description of Auslander}}

We have $M=\bigoplus_{i=0}^{\ell-1}\tau_n^{-i}\Lambda$ by Proposition \ref{tau n}.
We have \[\Hom_\Lambda(\tau_n^{-i}\Lambda,\tau_n^{-j}\Lambda)\stackrel{{\rm Prop.}\ \ref{tau n}}{\simeq}\left\{\begin{array}{cc}
\Hom_\Lambda(\tau_n^{-(i-j)}\Lambda,\Lambda)\stackrel{{\rm Lem.}\ \ref{tau and S}{\rm(b)}}{=}0&i>j,\\
\Hom_\Lambda(\Lambda,\tau_n^{-(j-i)}\Lambda)
\simeq\tau_n^{j-i}\Lambda&i\le j.
\end{array}\right.\]
Since we have
$\tau_n^{-i}\Lambda\simeq\SSS_n^{-i}\Lambda\simeq T^{\otimes_{\Lambda}i}$
for $0\le i<\ell$ by Lemma \ref{tau and S}(b) and \cite[2.12]{IO1},
we have the assertion.
\qed

\subsection{Proof of results in Section \ref{section2}}

Our Proposition \ref{homogeneous} is an immediate consequence of the following observation.

\begin{lemma}
Assume $\Hom_\Lambda(P_i,P_j)\neq0$, or equivalently $\Hom_\Lambda(I_i,I_j)\neq0$.
Then we have $\ell_i\le\ell_j$ and $\ell_{\sigma^{-1}(i)}\ge\ell_{\sigma^{-1}(j)}$.
\end{lemma}

\begin{proof}
If $\ell_i>\ell_j$, then by Theorem \ref{tau n}, we have
\[0\neq\Hom_\Lambda(I_i,I_j)\simeq\Hom_\Lambda(\tau_n^{\ell_j-1}I_i,P_{\sigma(j)})\]
by applying $\tau_n^{\ell_j-1}$.
By Lemma \ref{tau and S}(b), the right hand side is zero, a contradiction.
Thus we have $\ell_i\le\ell_j$.

If $\ell_{\sigma^{-1}(i)}<\ell_{\sigma^{-1}(j)}$, then we have
\[0\neq\Hom_\Lambda(P_i,P_j)\simeq\Hom_\Lambda(I_{\sigma^{-1}(i)},\tau_n^{-(\ell_{\sigma^{-1}(i)}-1)}P_j)\]
by applying $\tau_n^{-(\ell_{\sigma^{-1}(i)}-1)}$.
By Lemma \ref{tau and S}(b), the right hand side is zero, a contradiction.
\end{proof}

\medskip
Next we give a proof of Proposition \ref{psi and sigma}.

We have
\begin{equation}\label{2.2.1}
D\Hom_{\DD}(\SSS_n^jP_i,\Lambda)\simeq
\Hom_{\DD}(\Lambda,\SSS\SSS_n^jP_i)\simeq\Hom_{\DD}(\Lambda,\SSS_n^jI_i)
\simeq\Hom_{\DD}(\Lambda,\SSS_n^{j-\ell_i}P_{\sigma(i)}).
\end{equation}
Since we have functorial isomorphisms
\begin{eqnarray}\label{2.2.2}
\Pi\otimes_\Lambda-\simeq\bigoplus_{j\in\Z}\Hom_{\DD}(\Lambda,\SSS_n^j-)&:&\mod\Lambda\to\mod\Pi,\\ \label{2.2.3}
\Hom_\Pi(\Pi\otimes_\Lambda-,\Pi)\simeq\bigoplus_{j\in\Z}\Hom_{\DD}(\SSS_n^j-,\Lambda)&:&\mod\Lambda\to\mod\Pi^{\rm op}
\end{eqnarray}
on $\add\Lambda$, we have
\[\psi^*(\Pi\otimes_\Lambda P_i)\stackrel{\eqref{2.2.3}}{\simeq}\bigoplus_{j\in\Z}D\Hom_{\DD}(\SSS_n^jP_i,\Lambda)
\stackrel{\eqref{2.2.1}}{\simeq}\bigoplus_{j\in\Z}\Hom_{\DD}(\Lambda,\SSS_n^{j-\ell_i}P_{\sigma(i)})
\stackrel{\eqref{2.2.2}}{\simeq}\Pi\otimes_\Lambda P_{\sigma(i)}.\]
\qed

\medskip
Now we give a proof of Proposition \ref{phi is Nakayama}.

We always identify $\phi^*\Lambda$ with $\Lambda$ by the isomorphism
$\phi^*\Lambda=\Lambda_\phi\otimes_\Lambda\Lambda\ni \lambda\otimes\gamma\mapsto\lambda\cdot\phi\gamma\in\Lambda$.

(i) First we define a $K$-linear automorphism $\widetilde{\phi}$ of $\Pi$ by using the action of the autofunctor $\phi^*$ of $\DD$ on $\Hom_{\DD}(\Lambda,\SSS_n^{-i}\Lambda)$.

By Lemma \ref{tau and S}(c), we have an isomorphism $\rho:\phi^*\SSS_n^{-1}\to\SSS_n^{-1}\phi^*$ of autofunctors of $\DD$.
For any $i\ge0$, we define an isomorphism $\rho^i:\phi^*\SSS_n^{-i}\to\SSS_n^{-i}\phi^*$ of autofunctors of $\DD$ by $\rho^0=\id_{\phi^*}$ and
\[\rho^i:=(\phi^*\SSS_n^{-i}\xrightarrow{\rho_{\SSS_n^{1-i}}}\SSS_n^{-1}\phi^*\SSS_n^{1-i}\xrightarrow{\SSS_n^{-1}\rho_{\SSS_n^{2-i}}}\cdots
\xrightarrow{\SSS_n^{2-i}\rho_{\SSS_n^{-1}}}\SSS_n^{1-i}\phi^*\SSS_n^{-1}\xrightarrow{\SSS_n^{1-i}\rho}\SSS_n^{-i}\phi^*).\]
Clearly we have
\begin{equation}\label{i+j}
\rho^{i+j}=(\phi^*\SSS_n^{-i-j}\xrightarrow{\rho^i_{\SSS_n^{-j}}}\SSS_n^{-i}\phi^*\SSS_n^{-j}\xrightarrow{\SSS_n^{-i}\rho^j}\SSS_n^{-i-j}\phi^*)
\end{equation}
for any $i,j\ge0$.
For any $f\in\Hom_{\DD}(\Lambda,\SSS_n^{-i}\Lambda)$, define $\widetilde{\phi}f\in\Hom_{\DD}(\Lambda,\SSS_n^{-i}\Lambda)$ by
\[\widetilde{\phi}f:=(\Lambda=\phi^*\Lambda\xrightarrow{\phi^*f}\phi^*\SSS_n^{-i}\Lambda\xrightarrow{\rho^i_{\Lambda}}\SSS_n^{-i}\phi^*\Lambda=\SSS_n^{-i}\Lambda).\]
Clearly the restriction of $\widetilde{\phi}$ to $\Lambda$ coincides with $\phi$.

(ii) Next we show that $\widetilde{\phi}$ is a $K$-algebra automorphism of $\Pi$.

Fix $f\in\Hom_{\DD}(\Lambda,\SSS_n^{-i}\Lambda)$ and $g\in\Hom_{\DD}(\Lambda,\SSS_n^{-j}\Lambda)$.
Since $\rho^i$ is an isomorphism of functors, we have a commutative diagram
\begin{equation}\label{square}
\xymatrix{
\phi^*\SSS_n^{-i}\Lambda\ar^{\phi^*\SSS_n^{-i}g}[rr]\ar^{\rho^i_{\Lambda}}[d]
&&\phi^*\SSS_n^{-i-j}\Lambda\ar^{\rho^i_{\SSS_n^{-j}\Lambda}}[d]\\
\SSS_n^{-i}\phi^*\Lambda\ar^{\SSS_n^{-i}\phi^*g}[rr]&&\SSS_n^{-i-j}\phi^*\Lambda.
}\end{equation}
Now we have
\begin{eqnarray*}
\widetilde{\phi}f\cdot\widetilde{\phi}g&=&(\Lambda\xrightarrow{\widetilde{\phi}f}\SSS_n^{-i}\Lambda\xrightarrow{\SSS_n^{-i}(\widetilde{\phi}g)}\SSS_n^{-i-j}\Lambda)\\
&=&(\Lambda=\phi^*\Lambda\xrightarrow{\phi^*f}\phi^*\SSS_n^{-i}\Lambda
\xrightarrow{\rho^i_{\Lambda}}\SSS_n^{-i}\phi^*\Lambda
\xrightarrow{\SSS_n^{-i}\phi^*g}\SSS_n^{-i}\phi^*\SSS_n^{-j} \Lambda =\\&& \SSS_n^{-i}\phi^*\SSS_n^{-j} \Lambda
\xrightarrow{\SSS_n^{-i}\rho^j_{\Lambda}}\SSS_n^{-i-j}\phi^*\Lambda=\SSS_n^{-i-j}\Lambda)\\
&\stackrel{\eqref{square}}{=}&
(\Lambda=\phi^*\Lambda\xrightarrow{\phi^*f}\phi^*\SSS_n^{-i}\Lambda 
\xrightarrow{\phi^*\SSS_n^{-i}g}\phi^*\SSS_n^{-i-j}\Lambda
\xrightarrow{\rho^i_{\SSS_n^{-j}\Lambda}}\SSS_n^{-i}\phi^*\SSS_n^{-j} \Lambda =\\&& \SSS_n^{-i}\phi^*\SSS_n^{-j} \Lambda
\xrightarrow{\SSS_n^{-i}\rho^j_{\Lambda}}\SSS_n^{-i-j}\phi^*\Lambda=\SSS_n^{-i-j}\Lambda)\\
&\stackrel{\eqref{i+j}}{=}&
(\Lambda=\phi^*\Lambda\xrightarrow{\phi^*(f\SSS_n^{-i}g)}\phi^*\SSS_n^{-i-j}\Lambda
\xrightarrow{\rho^{i+j}_{\Lambda}}\SSS_n^{-i-j}\phi^*\Lambda=\SSS_n^{-i-j}\Lambda)\\
&=&\widetilde{\phi}(f\SSS_n^{-i}g).
\end{eqnarray*}

(iii) Since $\SSS^\ell\simeq[n(\ell-1)]\circ\phi^*$, we have
\begin{eqnarray}\label{SS phi}
\SSS\circ\SSS_n^{\ell-1}\simeq\phi^*.
\end{eqnarray}
By Serre duality, we have isomorphisms
\begin{equation}\label{by serre}
D\Hom_{\DD}(\Lambda,\SSS_n^{-i}\Lambda)\simeq\Hom_{\DD}(\SSS_n^{-i}\Lambda,\SSS\Lambda)
\simeq\Hom_{\DD}(\Lambda,\SSS\SSS_n^i\Lambda)
\stackrel{\eqref{SS phi}}{\simeq}\Hom_{\DD}(\Lambda,\phi^*\SSS_n^{i+1-\ell}\Lambda).
\end{equation}
Thus we have an isomorphism
\begin{equation}\label{DPi and Pi}
D\Pi=\bigoplus_{i\in\Z}D\Hom_{\DD}(\Lambda,\SSS_n^{-i}\Lambda)
\stackrel{\eqref{by serre}}{\simeq}\bigoplus_{i\in\Z}\Hom_{\DD}(\Lambda,\phi^*\SSS_n^{-i}\Lambda)=:N
\end{equation}
of $\Pi\otimes\Pi^{\rm op}$-modules,
where $\Pi$ acts on $N$ from the left in a usual way, and from the right by
\begin{equation}\label{right structure}
f\cdot g=(\Lambda\xrightarrow{f}\phi^*\SSS_n^{-i}\Lambda\xrightarrow{\phi^*\SSS_n^{-i}g}\phi^*\SSS_n^{-i-j}\Lambda)
\end{equation}
for any $f\in\Hom_{\DD}(\Lambda,\phi^*\SSS_n^{-i}\Lambda)\subset N$ and $g\in\Hom_{\DD}(\Lambda,\SSS_n^{-j}\Lambda)\subset\Pi$.
Let
\begin{equation}\label{bimodule iso}
\gamma:N=\bigoplus_{i\in\Z}\Hom_{\DD}(\Lambda,\phi^*\SSS_n^{-i}\Lambda)\to\Pi_{\widetilde{\phi}}=\bigoplus_{i\in\Z}\Hom_{\DD}(\Lambda,\SSS_n^{-i}\Lambda),
\end{equation}
be a $K$-linear isomorphism given by
\begin{equation*}
\gamma f:=(\Lambda\xrightarrow{f}\phi^*\SSS_n^{-i}\Lambda\xrightarrow{\rho^i_\Lambda}\SSS_n^{-i}\phi^*\Lambda=\SSS_n^{-i}\Lambda)
\end{equation*}
for any $f\in\Hom_{\DD}(\Lambda,\phi^*\SSS_n^{-i}\Lambda)$.
Then $\gamma$ is an isomorphism of $\Pi\otimes\Pi^{\rm op}$-modules 
by the following equalities for any $f\in\Hom_{\DD}(\Lambda,\phi^*\SSS_n^{-i}\Lambda)\subset N$ and $g\in\Hom_{\DD}(\Lambda,\SSS_n^{-j}\Lambda)\subset\Pi$.
\begin{eqnarray*}
\gamma f\cdot \widetilde{\phi} g
&=&(\Lambda\xrightarrow{f}\phi^*\SSS_n^{-i}\Lambda\xrightarrow{\rho^i_\Lambda}\SSS_n^{-i}\phi^*\Lambda
\xrightarrow{\SSS_n^{-i}\phi^*g}\SSS_n^{-i}\phi^*\SSS_n^{-j}\Lambda
\xrightarrow{\SSS_n^{-i}\rho^j_\Lambda}\SSS_n^{-i-j}\phi^*\Lambda=\SSS_n^{-i-j}\Lambda)\\
&\stackrel{\eqref{square}}{=}&
(\Lambda\xrightarrow{f}\phi^*\SSS_n^{-i}\Lambda\xrightarrow{\phi^*\SSS_n^{-i}g}\phi^*\SSS_n^{-i-j}\Lambda
\xrightarrow{\rho^i_{\SSS_n^{-j}\Lambda}}\SSS_n^{-i}\phi^*\SSS_n^{-j}\Lambda
\xrightarrow{\SSS_n^{-i}\rho^j_\Lambda}\SSS_n^{-i-j}\phi^*\Lambda=\SSS_n^{-i-j}\Lambda)\\
&\stackrel{\eqref{i+j}}{=}&(\Lambda\xrightarrow{f}\phi^*\SSS_n^{-i}\Lambda\xrightarrow{\phi^*\SSS_n^{-i}g}\phi^*\SSS_n^{-i-j}\Lambda
\xrightarrow{\rho^{i+j}_\Lambda}\SSS_n^{-i-j}\phi^*\Lambda=\SSS_n^{-i-j}\Lambda)\\
&\stackrel{\eqref{right structure}}{=}&\gamma(f\cdot g).
\end{eqnarray*}
Combining \eqref{DPi and Pi} and \eqref{bimodule iso}, we have an isomorphism $D\Pi\simeq\Pi_{\widetilde{\phi}}$ of $\Pi\otimes\Pi^{\rm op}$-modules.
Thus we have that $\widetilde{\phi}$ is a Nakayama automorphism of $\Pi$.
\qed

\medskip
Finally we give a proof of Theorem \ref{characterization}.

Assume that $\Lambda$ is $\ell$-homogeneous.
By Proposition \ref{phi is Nakayama}, the $K$-automorphism $\phi$ of $\Lambda$
extends to a Nakayama automorphism $\widetilde{\phi}$ of $\Pi$,
which satisfies $\widetilde{\phi}(I)=I$. Since $\widetilde{\phi}$ and $\psi$ coincide
in the outer automorphism group of $\Pi$, we have $\psi(I)=I$.

Assume $\psi(I)=I$.
Then $\psi$ induces a $K$-algebra automorphism $\psi|_\Lambda$ of $\Lambda$.
By Proposition \ref{psi and sigma}, we have $\psi|_\Lambda^*(P_i)\simeq P_{\sigma(i)}$ and $\psi|_\Lambda^*(I_i)\simeq I_{\sigma(i)}$ by Lemma \ref{tau and S}(c).
Since we have
$\psi|_\Lambda^*(\tau_n^jI_i)=\tau_n^j\psi|_\Lambda^*(I_i)=\tau_n^jI_{\sigma(i)}$ for any $j$ by Lemma \ref{tau and S}(c),
we have $\ell_i=\ell_{\sigma(i)}$.
By Proposition \ref{homogeneous}, we have that $\Lambda$ is homogeneous.
\qed


\begin{thebibliography}{15}
\bibitem[Am]{Am} C. Amiot, \emph{Cluster categories for algebras of global dimension 2 and quivers with potential},
arXiv:0805.1035, to appear in Ann. Inst. Fourier.
\bibitem[Au]{Au} M. Auslander, \textit{Representation dimension of Artin algebras}, Lecture notes, Queen Mary College, London, 1971.
\bibitem[Ba]{Ba} \textit{Spectral Methods in Representation Theory of Algebras and Applications to the Study of Rings of Singularities},
September 7-12, 2008, Banff International Research Station, Canada.
\bibitem[Bi]{Bi} \textit{The ADE chain}, October 31-November 1, 2008, Universitaet Bielefeld, Germany.
\bibitem[BK]{BK} A. I. Bondal, M. M. Kapranov, \textit{Representable functors, Serre functors, and reconstructions}, Izv. Akad. Nauk SSSR Ser. Mat. 53 (1989),  no. 6, 1183--1205, 1337;  translation in  Math. USSR-Izv.  35 (1990),  no. 3, 519--541
\bibitem[BIRS]{BIRS} A. Buan, O. Iyama, I. Reiten, J. Scott, \textit{Cluster structure for 2-Calabi-Yau categories and unipotent groups}, Compos. Math. 145 (2009), no. 4, 1035--1079.
\bibitem[BMRRT]{BMRRT} A. Buan, R. Marsh, I. Reiten, M. Reineke, G. Todorov, \textit{Tilting theory and cluster combinatorics}, Adv. in Math. 204 (2006), no. 2, 572--618.
\bibitem[D]{D} J. A. de la Pena, lectures at conferences \cite{Ba,P}.
\bibitem[DL]{DL} H. Lenzing, J. A. de la Pena, \textit{Extended canonical algebras and Fuchsian singularities}, arXiv:math.RT/0611532.
\bibitem[EH]{EH} K. Erdmann, T. Holm, \textit{Maximal $n$-orthogonal modules for selfinjective algebras}, Proc. Amer. Math. Soc. 136 (2008) no. 9, 3069--3078.
\bibitem[G]{G} P. Gabriel, \textit{Auslander-Reiten sequences and representation-finite algebras},
Representation theory, I (Proc. Workshop, Carleton Univ., Ottawa, Ont., 1979), pp. 1--71, Lecture Notes in Math., 831, Springer, Berlin, 1980.
\bibitem[GLS1]{GLS1} C. Geiss, B. Leclerc, J. Schr\"oer, \textit{Rigid modules over preprojective algebras}, Invent. Math. 165 (2006), no. 3, 589--632. 
\bibitem[GLS2]{GLS2} C. Geiss, B. Leclerc, J. Schr\"oer, \textit{Cluster algebra structures and semicanoncial bases for unipotent groups}, arXiv:math/0703039.
\bibitem[H]{H} D. Happel, \textit{Triangulated categories in the representation theory of finite-dimensional algebras},
London Mathematical Society Lecture Note Series, 119. Cambridge University Press, Cambridge, 1988.
\bibitem[HZ1]{HZ1} Z. Huang, X. Zhang, \textit{The existence of maximal $n$-orthogonal subcategories}, J. Algebra 321 (2009), no. 10, 2829--2842.
\bibitem[HZ2]{HZ2} Z. Huang, X. Zhang, \textit{Higher Auslander algebras admitting trivial maximal orthogonal subcategories}, arXiv:0903.0761.
\bibitem[HZ3]{HZ3} Z. Huang, X. Zhang, \textit{Trivial maximal $1$-orthognal subcategories for Auslander's $1$-Gorenstein algebras}, arXiv:0903.0762.
\bibitem[IIKNS]{IIKNS} R. Inoue, O. Iyama, A. Kuniba, T. Nakanishi, J. Suzuki, \textit{Periodicities of T-systems and Y-systems}, to appear in Nagoya Math. J., arXiv:0812.0667.
\bibitem[I1]{I1} O. Iyama, \textit{Higher-dimensional Auslander-Reiten theory on maximal orthogonal subcategories}, Adv. Math. 210 (2007), no. 1, 22--50.
\bibitem[I2]{I2} O. Iyama,  \textit{Auslander correspondence}, Adv. Math. 210 (2007), no. 1, 51--82.
\bibitem[I3]{I3} O. Iyama, \textit{Auslander-Reiten theory revisited}, Trends in Representation Theory of Algebras and Related Topics, 349--398, European Mathematical Society, 2008.
\bibitem[I4]{I4} O. Iyama, \textit{Cluster tilting for higher Auslander algebras}, arXiv:0809.4897.
\bibitem[IKM]{IKM} O. Iyama, K. Kato, J. Miyachi, \textit{Recollement of homotopy categories and Cohen-Macaulay modules}, arXiv:0911.0172.
\bibitem[IO1]{IO1} O. Iyama, S. Oppermann, \textit{$n$-representation-finite algebras and $n$-APR tilting}, to appear in Trans. Amer. Math. Soc., arXiv:0909.0593.
\bibitem[IO2]{IO2} O. Iyama, S. Oppermann, \textit{Stable categories over higher preprojective algebras}, arXiv:0912.3412.
\bibitem[IR]{IR} O. Iyama, I. Reiten, \textit{Fomin-Zelevinsky mutation and tilting modules over Calabi-Yau algebras}, Amer. J. Math. 130 (2008), no. 4, 1087--1149.
\bibitem[KST]{KST} H. Kajiura, K. Saito, A. Takahashi, \emph{Matrix factorization and representations of quivers. II. Type $ADE$ case}, Adv. Math. 211 (2007), no. 1, 327--362.
\bibitem[Ke1]{Ke1} B. Keller, \textit{Calabi-Yau triangulated categories}, Trends in Representation Theory of Algebras and Related Topics, 467--489, European Mathematical Society, 2008.
\bibitem[Ke2]{Ke2} B. Keller, \textit{Cluster algebras, quiver representations
  and triangulated categories}, arXiv:0807.1960.
\bibitem[Ke3]{Ke3} B. Keller, \textit{Deformed Calabi-Yau completions}, arXiv:0908.3499.
\bibitem[KR]{KR} B. Keller, I. Reiten, \textit{Cluster-tilted algebras are Gorenstein and stably Calabi-Yau},  Adv. Math. 211 (2007), 123-151.
\bibitem[Ko]{Ko} M. Kontsevich, \textit{Triangulated categories and geometry}, Course at the \'Ecole Normale Sup\'erieure, Paris, 1998.
\bibitem[Lada]{Lada} M. Lada, \textit{Relative homology and maximal $l$-orthogonal modules}, J. Algebra 321 (2009), no. 10, 2798--2811.
\bibitem[Ladk]{Ladk} S. Ladkani, lectures at conferences \cite{Ba,Bi,P}.
\bibitem[Le]{Le} H. Lenzing, lectures at conferences \cite{Ba,Bi,T}.
\bibitem[LM]{LM} H. Lenzing, H. Meltzer, \textit{Sheaves on a weighted projective line of genus one, and representations of a tubular algebra},
Representations of algebras (Ottawa, ON, 1992), 313--337, CMS Conf. Proc., 14, Amer. Math. Soc., Providence, RI, 1993.
\bibitem[Mi]{Mi} H. Minamoto \textit{Ampleness of two-sided tilting complexes}, preprint.
\bibitem[MY]{MY} J. Miyachi, A. Yekutieli, \textit{Derived Picard groups of finite-dimensional hereditary algebras}, Compositio Math. 129 (2001), no. 3, 341--368. 
\bibitem[P]{P} \textit{Workshop on Triangulated Categories and Singularities}, May 25-30, 2009, Universitaet Paderborn, Germany.
\bibitem[R]{R} C. M. Ringel, \textit{Tame algebras and integral quadratic forms}, Lecture Notes in Mathematics, 1099. Springer-Verlag, Berlin, 1984.
\bibitem[T]{T} \textit{Homological and geometric methods in algebra}, August 10-14, 2009, Norwegian University of Technology and Science, Norway.
\bibitem[Y]{Y} K. Yanagawa, lectures at Nagoya University, 2009.
\end{thebibliography}
\end{document}